\documentclass[12pt]{amsart}
\usepackage{amsmath,amsfonts,amssymb,latexsym}
\usepackage{hyperref}
\usepackage{graphicx}
\usepackage[a4paper, left=2cm, right=2cm, top=2cm, bottom=2cm]{geometry}

\newcommand{\RR}{\mathbb{R}}
\newcommand{\CC}{\mathbb{C}}
\newcommand{\NN}{\mathbb{N}}
\newcommand{\ZZ}{\mathbb{Z}}

\newcommand{\EE}{\mathbb{E}}
\newcommand{\OO}{\mathcal{O}}
\newcommand{\Oo}{\mathcal{O}}

\newcommand{\Bo}{\widehat{\mathcal{B}}}
\newcommand{\RE}{ {\rm Re \,} }

\newtheorem{Tw}{Theorem}

\newtheorem{Stw}{Proposition}

\theoremstyle{definition}
\newtheorem{Df}{Definition}

\theoremstyle{remark}
\newtheorem{Uw}{Remark}
\newtheorem{Prz}{Example}

\begin{document}
\keywords{PDEs with variable coefficients, moment functions, Borel summability, Stokes phenomenon, hyperfunctions}
\subjclass[2010]{35C10, 35C20, 40G10}
\title[The Stokes phenomenon for PDEs]{The Stokes phenomenon for certain partial differential equations with variable coefficients depending on the time variable}

\author{Bo\.{z}ena Tkacz}
\address{Faculty of Mathematics and Natural Sciences,
College of Science\\
Cardinal Stefan Wyszy\'nski University\\
W\'oycickiego 1/3,
01-938 Warszawa, Poland}
\email{b.tkacz@uksw.edu.pl}

\begin{abstract}
We consider the Stokes phenomenon for the solutions of some partial differential equations with variable coefficients in two complex variables, where initial data are holomorphic. We use the theory of (moment) summability and the theory of hyperfunctions to describe Stokes lines and jumps across these Stokes lines.
\end{abstract}

\maketitle 

\begin{section}{Introduction}
In this paper,  the main focus is on the Cauchy problem for partial differential equation with variable coefficients in two complex variables $(t,z)$
$$
\left\{ \begin{array}{ll} \partial_{t}u(t,z)=a(\partial_t t)^p t^q \partial_z^r u(t,z),\\
u(0,z)=\varphi(z)
\end{array}\right.
$$
with $p,q,r\in\NN_0, a\in\CC$ and $\varphi(z)\in\OO(D).$ Our purpose is to investigate the Stokes phenomenon for the formal solution $\widehat u(t,z)$ of the above equation i.e. we are interested in finding Stokes lines and describing jumps across these Stokes lines in terms of Laplace type hyperfunctions. The most important cases among those we considered are $r = 0, p\geq2$ and $r> 1, p\geq0$ (see main results -- Theorem \ref{tw1} and Theorem \ref{tw2}). In each of them we apply $m$-moment Borel transform $\Bo_m$ to $\widehat u(t,z)$ and $m$-moment Laplace transform in a nonsingular direction $d$ to $\Bo_m\widehat u$, hence we give an integral representation of a multisum of $\widehat u(t,z)$ in the direction $d$. Then we conclude what Stokes lines are and we derive the formula of the jumps.

It is worth pointing out that this work is a continuation of the paper \cite{Mic-Tk} in which we study the Stokes phenomenon for general homogeneous linear moment partial differential equation with constant coefficients.

\end{section}

\begin{section}{Notation}
A set 
\begin{displaymath}
S=S_d(\alpha,R)=\{z\in\tilde{\CC}\colon\ z=r\*e^{i\phi},\ r\in(0,R),\ \phi\in(d-\alpha/2,d+\alpha/2)\}.
\end{displaymath}
is called a~\emph{sector $S$ in a direction $d\in\RR$ with an opening $\alpha>0$ and a radius $R\in\RR_+$} in the universal covering space $\tilde{\CC}$ of
$\CC\setminus\{0\}$.\\
For $R=+\infty$, the sector $S$ is called \emph{unbounded} and one can write $S=S_d(\alpha)$ for short.\\
If the opening $\alpha$ is not essential, then for the unbounded sector $S$ the notation $S=S_d$ will be used.

A \emph{complex disc $D_r$ in $\CC$ with a radius $r>0$} is defined by
 \begin{displaymath}
 D_r=\{z\in\CC:|z|<r\}.
 \end{displaymath}
 If the radius $r$ is not essential, then one can abbreviate $D_r$ to $D$. 
 
A \emph{disc-sector} $S_d(\alpha)\cup D$ (resp. $S_d\cup D$) is designate briefly by $\widehat{S}_d(\alpha)$ (resp. $\widehat{S}_d$).

The notation $f\in\OO(G)$ means that a function $f$ is holomorphic on a domain $G\subset\CC^n$.
 
The set of all $\EE$-valued holomorphic functions on a domain $G\subseteq\CC^n$ will be represented by the symbol $\Oo(G,\EE)$, where $\EE$ is a complex Banach space with a norm $\|\cdot\|_{\EE}$. \\
 For simplicity, if $\EE=\CC$, then the set  $\Oo(G,\CC)$ will be written as $\Oo(G)$.

Through this paper, $\EE[[t]]$ stands for the space of all formal power series $\hat{u}=\sum_{n=0}^{\infty} a_{n} t^{n}$ with $a_n\in\EE$.\\
If the formal power series $\hat{u}$ is convergent, then its sum will be denoted by $u$.

\begin{Df}
Assume that $k\in\RR$, $S$ is an unbounded sector and $u\in\OO(S,\EE)$. The function $u$ is of \emph{exponential growth of order at most $k$},
if for every proper subsector $S^*\prec S$ (i.e. $\overline{S^{*}}\setminus\{0\} \subseteq S$) there exist constants $C_1, C_2>0$ such that
$\|u(x)\|_{\EE}\le C_1\*e^{C_2|x|^{k}}$ for every $x\in S^*$. If this is so, then one can write $u\in\OO^{k}(S,\EE)$.
\end{Df}
\end{section}

\begin{section}{Kernel and moment functions, k-summability}
In this section we review some of the standard definitions and facts on moment methods introduced by Balser \cite{B2}. Also as in \cite{Mic-Tk} we give definitions of Gevrey order, moment Borel transforms and  $k$-summability.

\begin{Df}[see {\cite[Section 5.5]{B2}}]
    \label{df:moment}
    A pair of functions $e_m$ and $E_m$ is called \emph{kernel functions of order $k$} ($k>1/2$) if
     the following properties are satisfied
   \begin{enumerate}
    \item[1.] $e_m\in\Oo(S_0(\pi/k))$.
    \item[2.] $e_m(z)/z$ is integrable at the origin.
    \item[3.] $e_m(x)\in\RR_+$ for $x\in\RR_+$.
    \item[4.] $e_m$ is exponentially flat of order $k$ as $z\to\infty$ in $S_0(\pi/k)$ (i.e. for every $\varepsilon > 0$ there exist ${A,B > 0}$
     such that $|e_m(z)|\leq A e^{-(|z|/B)^k}$ for $z\in S_0(\pi/k-\varepsilon)$).
    \item[5.] $E_m\in\Oo^{k}(\CC)$.
    \item[6.] $E_m(1/z)/z$ is integrable at the origin in $S_{\pi}(2\pi-\pi/k)$.
    \item[7.] \emph{The corresponding moment function $m$ of order $1/k$} is given by the Mellin transform of $e_m$
     \begin{gather}
      \label{eq:e_m}
      m(u):=\int_0^{\infty}x^{u-1} e_m(x)dx \quad \textrm{for} \quad \RE u \geq 0,
     \end{gather}
      and the kernel function $E_m$ has the power series expansion
     \begin{gather}
      \label{eq:E_m}
      E_m(z)=\sum_{n=0}^{\infty}\frac{z^n}{m(n)} \quad  \textrm{for} \quad z\in\CC.
     \end{gather}
     
    So the function $m$ describes the connection between kernel functions $e_m$ and $E_m$. 
    \item[8.] $m(0)=1$.
   \end{enumerate}
   \end{Df}

\begin{Uw}
One can easily seen that the moment function $m$ uniquely determines the kernel functions $e_m$ and $E_m$.
\end{Uw}

Observe that if $k\leq 1/2$,  then the set $S_{\pi}(2\pi-\pi/k)$ is not defined. So in that case Definition \ref{df:moment} is not valid. 
Thus, the kernel functions of order $k\leq 1/2$ and the corresponding moment functions must be defined in another way.
To this end we use the ramification at $z=0$.

   \begin{Df}[see {\cite[Section 5.6]{B2}}]
     \label{df:small}
   If one can find a pair of kernel functions $e_{\widetilde{m}}$ and $E_{\widetilde{m}}$ of order $pk>1/2$ (for some $p\in\NN$) so that
     \begin{gather*}
      e_m(z)=e_{\widetilde{m}}(z^{1/p})/p \quad \textrm{for} \quad z\in S_0(\pi/k).
     \end{gather*}
     then a function $e_m$ is called a \emph{kernel function of order $k>0$}.
     \emph{The corresponding moment function $m$ of order $1/k>0$} for a given kernel function $e_m$ of order $k>0$ is defined by (\ref{eq:e_m}) and 
     the \emph{kernel function $E_m$ of order $k>0$} by (\ref{eq:E_m}).
    \end{Df}
    
    \begin{Uw}
    Notice that by Definitions \ref{df:moment} and \ref{df:small} we obtain
     \begin{eqnarray*}
      m(u)=\widetilde{m}(pu) & \textrm{and} &
      E_m(z)=\sum_{j=0}^{\infty}\frac{z^j}{m(j)}=\sum_{j=0}^{\infty}\frac{z^j}{\widetilde{m}(jp)}.
     \end{eqnarray*}
    \end{Uw}

Now, similarly as in \cite{Mic8}, we define a moment function with a real order.

\begin{Df}
 \label{df:moment_general}
     A function $m$ is a \emph{moment function of order $1/k<0$},  if $1/m$ is a moment function of order $-1/k>0$.
     \par
     A function $m$ is a \emph{moment function of order $0$},  if there exist moment functions $m_1$ and $m_2$ of the same order $1/k>0$ such that $m=m_1/m_2$.
\end{Df}

Based on Definition \ref{df:moment_general} and by \cite[Theorems 31 and 32]{B2} we get
\begin{Stw}
\label{propo:1}
 Assume that $m_1$, $m_2$ are moment functions of orders $k_1,k_2\in\RR$ respectively. Then
 \begin{itemize}
  \item $m_1m_2$ is a moment function of order $s_1+s_2$,
  \item $m_1/m_2$ is a moment function of order $s_1-s_2$.
 \end{itemize}
\end{Stw}

Now, we give two examples of moment functions, which will be extensively used in this paper.
\begin{Prz}
\label{ex:functions}
 Let $k>0$. Then the kernel functions and the corresponding moment function,
 satisfying Definition \ref{df:moment} or \ref{df:small}, are $e_m(z)=kz^{k}e^{-z^k}$, $m(u)=\Gamma(1+u/k)$ and $E_m(z)=\sum_{j=0}^{\infty}\frac{z^j}{\Gamma(1+j/k)}=:\mathbf{E}_{1/k}(z)$, where $\mathbf{E}_{1/k}$ is the Mittag-Leffler function of index $1/k$.
\end{Prz}

\begin{Prz}
Let $s\in\RR$. By $\Gamma_s$ we denote the function
\[
  \Gamma_s(u):=\left\{
  \begin{array}{lll}
    \Gamma(1+su) & \textrm{for} & s \geq 0\\
    1/\Gamma(1-su) & \textrm{for} & s < 0.
  \end{array}
  \right.
\]
Observe that the above function is a moment function of order $s\in\RR$.\\
 Moreover, every moment function $m$ of order $s$ has the same growth as $\Gamma_s$ (see {\cite[Section 5.5]{B2}}). Precisely speaking,
we have 
\end{Prz}

\begin{Stw}[see {\cite[Section 5.5]{B2}}]
  \label{pr:order}
  If $m$ is a moment function of order $s\in\RR$
  then there exist constants $a,A,c,C>0$ such that
    \begin{gather*}
    a c^n\Gamma_s(n)\leq m(n) \leq A C^n\Gamma_s(n) \quad \textrm{for every} \quad n\in\NN_0.
    \end{gather*}
  \end{Stw}

In the last part of this section we describe moment Borel transform, the Gevrey order and k-summability.

\begin{Df}
 Suppose that $m$ is a moment function. An \emph{$m$-moment Borel transform} is a linear operator $\Bo_{m}\colon \EE[[t]]\to\EE[[t]]$ defined by
 \[
  \Bo_{m}\bigg(\sum_{j=0}^{\infty}u_jt^{j}\bigg):=
  \sum_{j=0}^{\infty}\frac{u_j}{m(j)}t^{j}
 \]
\end{Df}

   \begin{Df}
    \label{df:summab}
    Let $s\in\RR$. A formal power series
    $\widehat{u}\in\EE[[t]]$ is called a \emph{formal power series of Gevrey order $s$}, if
    there exists a disc $D\subset\CC$ with centre at the origin such that
    $\Bo_{\Gamma_s}\widehat{u}\in\Oo(D,\EE)$. By the symbol $\EE[[t]]_s$ we denote the space of formal power series of Gevrey 
    order $s$.
   \end{Df}

\begin{Uw}
 Let $\widehat{u}\in\EE[[t]]_s$ and $s\leq 0$. Then the formal series $\widehat{u}$ is convergent, hence its sum $u$ is well defined.
 Moreover, $\widehat{u}\in\EE[[t]]_0 \Longleftrightarrow u\in\Oo(D,\EE)$ and
 $\widehat{u}\in\EE[[t]]_s \Longleftrightarrow u\in\Oo^{-1/s}(\CC,\EE)$ for $s<0$.
\end{Uw}

\begin{Df}
  Let $d\in\RR$. Assume also that $e_m, E_m$ is a pair of kernel functions of order $1/k>0$ with a moment function $m$.
  \begin{itemize}
  \item If $v\in\Oo^{k}(\widehat{S}_d,\EE)$ then the integral operator $T_{m,d}$ defined by
  \[
   (T_{m,d}v)(t):=\int_{e^{id}\RR_+}e_{m}\Big(\frac{s}{t}\Big)v(s)\frac{ds}{s}
  \]
 is called an \emph{$m$-moment Laplace transform in a direction $d$}.
 \item If $v\in\Oo(S_d(\frac{\pi}{k}+\varepsilon,R),\EE)$ for some $\varepsilon,R>0$ then the integral operator $T^-_{m,d}$ defined by
  \[
   (T^-_{m,d}v)(s):=-\frac{1}{2\pi i}\int_{\gamma(d)}E_{m}\Big(\frac{s}{t}\Big)v(t)\frac{dt}{t}
  \]
(where a path $\gamma(d)$ is the boundary of a sector contained in $S_d(\frac{\pi}{k}+\varepsilon,R)$ with bisecting direction $d$,
a finite radius, an opening slightly larger than $\pi/k$, and the orientation is negative) is called an \emph{inverse $m$-moment Laplace transform in a direction $d$}.
 \end{itemize}
\end{Df}

\begin{Uw}
Notice that $T_{m,d}(t^{n})=m(n)t^{n}$ for every $n\in\NN_0$, thus $ T_{m,d}\Bo_{m}u=u$
 for every $u\in\Oo(D)$.
\end{Uw}

\begin{Df}
Let $k>0$ and $d\in\RR$. The formal power series $\widehat{u}\in\EE[[t]]$ is called
\emph{$k$-summable in a direction $d$} if there exist $\varepsilon>0$ and a disc-sector $\widehat{S}_d=\widehat{S}_d(\varepsilon)$
in a direction $d$ such that
$v=\Bo_{\Gamma_{1/k}}\widehat{u}\in\Oo^k(\widehat{S}_d,\EE)$.
\par
 Moreover, the \emph{$k$-sum of $\widehat{u}$ in the direction $d$} is given by
\begin{equation}
\label{eq:sum}
       u^d(t)=\mathcal{S}_{k,d}\widehat{u}(t):=(T_{m,\theta}v)(t)=\int_{e^{i\theta}\RR_+}e_{m}\Big(\frac{s}{t}\Big)v(s)\frac{ds}{s}\quad\textrm{for}\quad
       \theta\in(d-\varepsilon/2,d+\varepsilon/2).
\end{equation}
\end{Df}

\begin{Df}
If $\widehat{u}\in\EE[[t]]$ is $k$-summable in all directions $d$ but (after identification modulo $2\pi$)
finitely many directions $d_1,\dots,d_n$ then
$\widehat{u}$ is called \emph{$k$-summable} and $d_1,\dots,d_n$ are called \emph{singular directions of $\widehat{u}$}.
\end{Df}
\end{section}

\begin{section}{The Stokes phenomenon via hyperfunctions}
In this section we recall the notion of the Stokes phenomenon for formal power series $\widehat{u}\in\EE[[t]]$ (see \cite[Definition 7]{Mic-Pod}).
\begin{Df}
\label{df:stokes}
Assume that $\widehat{u}\in\EE[[t]]$ is $k$-summable with finitely many singular directions $d_1,d_2,\dots,d_n$. Then for every $l=1,\dots,n$ a set $\mathcal{L}_{d_{l}}=\{t\in\tilde{\CC}\colon \arg t=d_{l}\}$  is called a \emph{Stokes line} for $\hat{u}$ and a set $\mathcal{L}_{d_l\pm\frac{\pi}{2k}}$ is called an \emph{anti-Stokes line} for $\hat{u}$.

Moreover, if $d_l^+$ (resp. $d_l^-$) denotes a direction close to $d_l$ and greater (resp. less) than $d_l$, and let $u^{d_l^+}:=\mathcal{S}_{k, d_l^+}\hat{u}$ (resp. $u^{d_l^-}:=\mathcal{S}_{k, d_l^-}\hat{f}$)
then the difference $J_{\mathcal{L}_{d_{l}}}\widehat{u}:=u^{d_l^+}- u^{d_l^-} $ is called a \emph{jump} for $\hat{u}$ across the Stokes line $\mathcal{L}_{d_l}$. 
\end{Df}

Now (see similar approach in \cite{Im, Mal3, S-S}) we describe jumps across the Stokes lines in terms of hyperfunctions.
For more information about the theory of hyperfunctions we refer the reader to \cite{Kaneko}.

As in \cite{Mic-Tk} we consider the space 
$$\mathcal{H}^k(\mathcal{L}_d):=\OO^k(D\cup (S_d\setminus \mathcal{L}_d))\Big/\OO^k(\widehat{S}_d)$$
of Laplace type hyperfunctions supported by $\mathcal{L}_d$ with exponential growth of order $k$. It means that every hyperfunction
$G\in\mathcal{H}^k(\mathcal{L}_d)$ may be written as
\[
G(s)=[g(s)]_{d}=\{g(s)+h(s)\colon h(s)\in\OO^{k}(\widehat{S}_d)\}
\]
for some defining function $g(s)\in\OO^k(D\cup (S_d\setminus \mathcal{L}_d))$. 

Let $\gamma_{d}$ be a path consisting of the half-lines from 
$e^{id^-}\infty$ to $0$ and from $0$ to $e^{id^+}\infty$, i.e.
$\gamma_{d}=-\gamma_{d^-}+\gamma_{d^+}$ with $\gamma_{d^{\pm}}=\mathcal{L}_{d^{\pm}}$.
By the K\"othe type theorem \cite{Kot} one can treat the hyperfunction $G(s)=[g(s)]_d$ as the analytic functional defined by
\begin{gather}
\label{eq:kothe}
G(s)[\varphi(s)]:=\int_{\gamma_d}g(s)\varphi(s)\,ds,
\end{gather}
for such small $\varphi\in\OO^{-k}(\widehat{S}_d)$
that the function $s\mapsto g(s)\varphi(s)$ belongs to the space
$\OO^{-k}(D\cup (S_d\setminus \mathcal{L}_d))$.

To describe the jumps across the Stokes lines in terms of hyperfunctions,
we assume that $\widehat{f}\in\CC[[t]]$ is $k$-summable, $m$ is a moment function of order $1/k$ and $d$ is a singular direction. 
By (\ref{eq:sum}) the jump for $\widehat{f}$ across the Stokes line $\mathcal{L}_d$ is given by 
\begin{equation}
\label{eq:jump_hyp}
 J_{\mathcal{L}_d}\widehat{f}(t)=f^{d^+}(t)-f^{d^-}(t)=(T_{m,d^+}-T_{m,d^-})\Bo_{m}\widehat{f}(t).
\end{equation}

Notice that we can treat
$g_0(t):=\Bo_{m}\widehat{f}(t)\in\Oo^{k}(D\cup (S_{d}\setminus \mathcal{L}_{d}))$
as a defining function of the hyperfunction $G_0(s):=[g_0(s)]_{d}\in \mathcal{H}^{k}(\mathcal{L}_{d})$.
Thus, by (\ref{eq:kothe}) with (\ref{eq:jump_hyp}) we have that 
\begin{gather}
\label{eq:J}
 J_{\mathcal{L}_d}\widehat{f}(t)=G_0(s)\Big[\frac{e_{m}(s/t)}{s}\Big]\quad \textrm{for sufficiently small}\ r>0\ \textrm{and}\
 t\in S_{d}\big(\frac{\pi}{k}, r\big).
\end{gather}

If $\mathcal{L}_d$ is a Stokes line for $k$-summable $\widehat{u}=\widehat{u}(t,z)\in\Oo(D)[[t]]$, then
we are able to describe jumps for $\widehat{u}(t,z)$ at the point $z=0$ in terms of hyperfunctions. Namely we have 
\begin{gather*}
J_{\mathcal{L}_{d}}\widehat{u}(t,0)=(T_{m,d}F_{0})(t)= F_{0}(s)\Big[\frac{e_{m}(s/t)}{s}\Big],\quad \textrm{where}\quad
F_{0}(s)=[\Bo_{m}\widehat{u}(s,0)]_{d}\in \mathcal{H}^k(\mathcal{L}_{d}).
\end{gather*}
\end{section}

\begin{Uw}
\label{re:z}
Since the derivation operator $\partial_z$ is invariant under the translation we are also able to describe jumps for $\widehat{u}(t,z)$ at any point $z\in D$ (see \cite[Remark 9]{Mic-Tk}, for more details).
 \end{Uw}
 
 \begin{section}{Partial differential equations with variable coefficients}
 In this section we will consider certain PDEs with variable coefficients -- we derive Stokes lines and jumps across these Stokes lines in terms of hyperfunctions.

Let us study the following equation
\begin{equation}
\label{eq}
\left\{ \begin{array}{ll} \partial_{t}u(t,z)=a(\partial_t t)^p t^q \partial_z^r u(t,z),\,\, \mathrm{with}\,\, p,q,r\in\NN_0, a\in\CC\\ 
u(0,z)=\varphi(z)\in\OO(D).
\end{array}\right.
\end{equation}

The above Cauchy problem has a unique formal solution $$\widehat{u}(t,z)=\sum_{n=0}^{\infty}a^n (n!)^{p-1}\,(q+1)^{n(p-1)}\,\varphi^{(nr)}(z)\,t^{n(q+1)}.$$

The cases $r=0$ and $r>0$ will be the ones of interest to us.\\

\textbf{Case 1.} Let $r=0$ and $m_1(n)=\Gamma\big(1+\frac{n}{q+1}\big)$ be a moment function of order $\frac{1}{q+1}$ corresponding to a kernel function $e_{m_1}(z)$ of order $q+1$. 

\begin{enumerate}
\renewcommand{\labelenumi}{(\alph{enumi})}
\item Assume that $p=0$. Then $$\widehat{u}(t,z)=\sum_{n=0}^{\infty}a^n (n!)^{-1}\,(q+1)^{-n}\,\varphi(z)\,t^{n(q+1)}=\varphi(z)\sum_{n=0}^{\infty}\frac{\big(a(q+1)^{-1}t^{q+1}\big)^n}{n!}=\varphi(z)e^{\frac{at^{q+1}}{q+1}}$$
is an entire function of $t$. Thus, there is no singular directions.

\item Suppose that $p=1$. Then  $$\widehat{u}(t,z)=\sum_{n=0}^{\infty}a^n \varphi(z)\,t^{n(q+1)}= \frac{\varphi(z)}{1-at^{q+1}},\,\,\,\,\mathrm{if}\,\,\, \big|at^{q+1}\big|<1$$
is an analytic function for $|t|<|a|^{-\frac{1}{q+1}}$, so there is no singular directions.

\item Let $p\geq2$. By Proposition \ref{propo:1} function $$m(n):=\underbrace{m_1(n)\cdot\ldots\cdot m_1(n)}_{(p-1)-\mathrm{times}}=\bigg(\Gamma\big(1+\frac{n}{q+1}\big)\bigg)^{p-1}$$ is a moment function of order $\frac{p-1}{q+1}$ corresponding to a kernel function $e_{m}(z)$ of order $\frac{q+1}{p-1}$.

First, we apply the $m$-moment Borel transform to $\widehat{u}(t,z)$. 
\begin{multline*}
(\Bo_m\widehat{u})(t,z)=\sum_{n=0}^{\infty}\frac{a^n (n!)^{p-1}\,(q+1)^{n(p-1)}\,\varphi(z)}{m(n(q+1))}\,t^{n(q+1)}= \sum_{n=0}^{\infty}\frac{a^n (n!)^{p-1}\,(q+1)^{n(p-1)}\,\varphi(z)}{\bigg(\Gamma\big(1+\frac{n(q+1)}{q+1}\big)\bigg)^{p-1}}\,t^{n(q+1)}\\
= \sum_{n=0}^{\infty}\frac{a^n (n!)^{p-1}\,(q+1)^{n(p-1)}\,\varphi(z)}{\bigg(\Gamma\big(1+n \big)\bigg)^{p-1}}\,t^{n(q+1)}=\sum_{n=0}^{\infty}\frac{a^n (n!)^{p-1}\,(q+1)^{n(p-1)}\,\varphi(z)}{(n!)^{p-1}}\,t^{n(q+1)}\\
=\frac{\varphi(z)}{1-a\,(q+1)^{p-1}\,t^{q+1}},
\end{multline*}
if $\big|a\,(q+1)^{p-1}\,t^{q+1}\big|<1.$

Let $f(s,z):=(\Bo_m\widehat{u})(s,z)$, then by using $m$-moment Laplace transform in a nonsingular direction $d$ we get 
$$(T_{m,d}f)(t,z)=\int_{e^{id}\RR_+}e_m\Big(\frac{s}{t}\Big)\*f(s,z)\frac{ds}{s}=\int_{e^{id}\RR_+}e_m\Big(\frac{s}{t}\Big)\*\*\frac{\varphi(z)}{1-a\,(q+1)^{p-1}\,s^{q+1}}\frac{ds}{s}.$$

Thus, if $d\neq\frac{2k\pi-\arg{a}}{q+1}$ for $k\in\ZZ$, then the unique formal solution $\widehat{u}(t,z)$ of this Cauchy problem is $\frac{q+1}{p-1}$-summable in the direction $d$ and for every $\varepsilon>0$ there exists $\tilde r>0$ such that
its $\frac{q+1}{p-1}$-sum $u\in\OO(S_{d}(\frac{\pi(p-1)}{q+1}-{\varepsilon},\tilde r)\times D)$ is given by 
$$
u(t,z)=u^{d}(t,z)=\int_{e^{id}\RR_+}e_m\Big(\frac{s}{t}\Big)\*\*\frac{\varphi(z)}{1-a\,(q+1)^{p-1}\,s^{q+1}}\frac{ds}{s}.
$$

Let $d_k:=\frac{2k\pi-\arg{a}}{q+1}$ for $k=0,\ldots,q$.
Then  $\mathcal{L}_{d_k}$ are Stokes lines for $\widehat{u}$ and the jump is given by 
\begin{multline*}
J_{\mathcal{L}_{d_k}}\widehat u(t,z)= u^{d_k+\varepsilon}(t,z) - u^{d_k-\varepsilon}(t,z)=F_z(s)\bigg[\frac{e_m(s/t)}{s}\bigg]=\big[\Bo_m\widehat{u}(s,z)\big]\bigg[\frac{e_m(s/t)}{s}\bigg]\\
=\bigg[\frac{\varphi(z)}{1-a\,(q+1)^{p-1}\,s^{q+1}}\bigg]_{d_k}\bigg[\frac{e_m(s/t)}{s}\bigg]\\
=\Bigg[\frac{\varphi(z)}{\big(1-a^{\frac{1}{q+1}}(q+1)^{\frac{p-1}{q+1}}s\big)\cdot\big(1-a^{\frac{1}{q+1}}(q+1)^{\frac{p-1}{q+1}}e^{\frac{2\pi i}{q+1}}s\big)\cdot\ldots\cdot\big(1-a^{\frac{1}{q+1}}(q+1)^{\frac{p-1}{q+1}}e^{\frac{2\pi iq}{q+1}}s\big)}\Bigg]_{d_k}\bigg[\frac{e_m(s/t)}{s}\bigg]\\
=2\pi i \frac{\varphi(z)}{\prod_{\substack{j=0\\ j\neq k}}^q \big(1-a^{\frac{1}{q+1}}(q+1)^{\frac{p-1}{q+1}}e^{\frac{2\pi ij}{q+1}}\big(\frac{1}{a}\big)^{\frac{1}{q+1}}(q+1)^{-\frac{p-1}{q+1}}\,e^{-\frac{2\pi ik}{q+1}}\big)}\frac{e_m\Big(\big(\frac{1}{a}\big)^{\frac{1}{q+1}}(q+1)^{-\frac{p-1}{q+1}}\,e^{-\frac{2\pi ik}{q+1}}t^{-1}\Big)}{\big(\frac{1}{a}\big)^{\frac{1}{q+1}}(q+1)^{-\frac{p-1}{q+1}}\,e^{-\frac{2\pi ik}{q+1}}}\\
= 2\pi i \frac{\varphi(z)}{\prod_{\substack{j=0\\ j\neq k}}^q \Big(1-(e^{2\pi i})^{\frac{j-k}{q+1}}\Big)} \frac{e_m\Big(\big(\frac{1}{a}\big)^{\frac{1}{q+1}}(q+1)^{-\frac{p-1}{q+1}}\,e^{-\frac{2\pi ik}{q+1}}t^{-1}\Big)}{\big(\frac{1}{a}\big)^{\frac{1}{q+1}}(q+1)^{-\frac{p-1}{q+1}}\,e^{-\frac{2\pi ik}{q+1}}}\\=
2\pi i \frac{\varphi(z)}{q+1} \frac{e_m\Big(\big(\frac{1}{a}\big)^{\frac{1}{q+1}}(q+1)^{-\frac{p-1}{q+1}}\,e^{-\frac{2\pi ik}{q+1}}t^{-1}\Big)}{\big(\frac{1}{a}\big)^{\frac{1}{q+1}}(q+1)^{-\frac{p-1}{q+1}}\,e^{-\frac{2\pi ik}{q+1}}}.
\end{multline*}

Observe that from \cite[Theorem 31]{B2} one can derive the function $e_m(z)$ by induction on $p$ ie.

\begin{itemize}
\item Set $p=2$, then $m^{\mathbf{0}}:=m(n)=m_1(n)$ and
$$e_{m^{\mathbf{0}}}(z)=e_{m_1}(z)\stackrel{\mathrm{by\,\,\, Example} \,\ref{ex:functions}}{=}(q+1)z^{q+1}e^{-z^{q+1}}.$$
 
\item Set $p=3$, then $m^{\mathbf{1}}(n):=m(n)=m_1(n)m^{\mathbf{0}}(n)$ and 
\begin{multline*}
 e_{m^{\mathbf{1}}}(z)=T_{m_1,d}\,\bigg(e_{m^{\mathbf{0}}}\Big(\frac{1}{u}\Big)\bigg)\Big(\frac{1}{z}\Big) =\int_{e^{id}\RR_+} e_{m_1}(uz)\,\,e_{m^{\mathbf{0}}}\Big(\frac{1}{u}\Big)\frac{du}{u}\\
 \stackrel{\mathrm{by\,\,\, Example}\, \ref{ex:functions}}{=}\int_{e^{id}\RR_+}(q+1)(uz)^{q+1}e^{-(uz)^{q+1}}(q+1)\Big(\frac{1}{u}\Big)^{q+1}e^{-\big(\frac{1}{u}\big)^{q+1}}\frac{du}{u}\\
 =\int_{e^{id}\RR_+}(q+1)^{2}z^{q+1}e^{-(uz)^{q+1}}e^{-\big(\frac{1}{u}\big)^{q+1}}\frac{du}{u}.
\end{multline*}

\item Set $p=4$, then $m^{\mathbf{2}}(n):=m(n)=m_1(n)m^{\mathbf{1}}(n)$ and 
\begin{multline*}
 e_{m^{\mathbf{2}}}(z)=T_{m_1,d}\,\bigg(e_{m^{\mathbf{1}}}\Big(\frac{1}{u}\Big)\bigg)\Big(\frac{1}{z}\Big) =\int_{e^{id}\RR_+} e_{m_1}(uz)\,\,e_{m^{\mathbf{1}}}\Big(\frac{1}{u}\Big)\frac{du}{u}\\
 \stackrel{\mathrm{by\,\,\, Example}\, \ref{ex:functions}}{=}\int_{e^{id}\RR_+}(q+1)(uz)^{q+1}e^{-(uz)^{q+1}}e_{m^{\mathbf{1}}}\Big(\frac{1}{u}\Big)\frac{du}{u}.
 \end{multline*}
 
 \item Hence, in general if $p\geq3$, then $m^{\mathbf{p-2}}(n):=m(n)=m_1(n)m^{\mathbf{p-3}}(n)$ and 
 \begin{multline*}
  e_{m^{\mathbf{p-2}}}(z)=T_{m_1,d}\,\bigg(e_{m^{\mathbf{p-3}}}\Big(\frac{1}{u}\Big)\bigg)\Big(\frac{1}{z}\Big) =\int_{e^{id}\RR_+} e_{m_1}(uz)\,\,e_{m^{\mathbf{p-3}}}\Big(\frac{1}{u}\Big)\frac{du}{u}\\
 \stackrel{\mathrm{by\,\,\, Example}\, \ref{ex:functions}}{=}\int_{e^{id}\RR_+}(q+1)(uz)^{q+1}e^{-(uz)^{q+1}}e_{m^{\mathbf{p-3}}}\Big(\frac{1}{u}\Big)\frac{du}{u}.
 \end{multline*}
\end{itemize}

\end{enumerate}
\bigskip

\textbf{Case 2.} Assume now that $r>0$. Let $m_1(n)=\Gamma\big(1+\frac{n}{q+1}\big)$ be a moment function of order $\frac{1}{q+1}$ corresponding to a kernel function $e_{m_1}(z)$ of order $q+1$, and $m_2(n)=\Gamma\big(1+\frac{nr}{q+1}\big)$ be a moment function of order $\frac{r}{q+1}$ corresponding to a kernel function $e_{m_2}(z)$ of order $\frac{q+1}{r}$.
\begin{enumerate}
\renewcommand{\labelenumi}{(\alph{enumi})}

\item Let $r=1$ and $p=0$. Then $$\widehat{u}(t,z)=\sum_{n=0}^{\infty}a^n (n!)^{-1}\,(q+1)^{-n}\,\varphi^{(n)}(z)\,t^{n(q+1)}=\sum_{n=0}^{\infty}\frac{\big(a(q+1)^{-1}t^{q+1}\big)^n\varphi^{(n)}(z)}{n!}
=\varphi\bigg(z+ \frac{at^{q+1}}{q+1}\bigg)$$
for sufficiently small $z$ and $t$ is an analytic function, so there is no singular directions.

\item Suppose that $r>1$, $p=0$ and $\varphi(z)\in\OO^{\frac{r}{r-1}}\Bigl(\widetilde{\CC\setminus\{z_0\}}\Bigr)$ for some $z_0\in\CC\setminus\{0\}$. Then by Proposition \ref{propo:1} function $m(n):=\frac{m_2(n)}{m_1(n)}=\frac{\Gamma\big(1+\frac{nr}{q+1}\big)}{\Gamma\big(1+\frac{n}{q+1}\big)}$  is a moment function of order $\frac{r-1}{q+1}$ corresponding to a kernel function $e_{m}(z)$ of order $\frac{q+1}{r-1}$.

We start with applying the $m$-moment Borel transform to $\widehat{u}(t,z)$. 
\begin{multline*}
(\Bo_m\widehat{u})(t,z)=\sum_{n=0}^{\infty}\frac{a^n (n!)^{-1}\,(q+1)^{-n}\,\varphi^{(nr)}(z)}{m(n(q+1))}\,t^{n(q+1)}\\
=\sum_{n=0}^{\infty}\frac{a^n (n!)^{-1}\,(q+1)^{-n}\,\varphi^{(nr)}(z)}{\Gamma\big(1+\frac{n(q+1)r}{q+1}\big)\Big(\Gamma\big(1+\frac{n(q+1)}{q+1}\big)\Big)^{-1}}\,t^{n(q+1)}=\sum_{n=0}^{\infty}\frac{a^n (n!)^{-1}\,(q+1)^{-n}\,\varphi^{(nr)}(z)}{\Gamma(1+nr)\Big(\Gamma(1+n)\Big)^{-1}}\,t^{n(q+1)}\\
=\sum_{n=0}^{\infty}\frac{a^n (n!)^{-1}\,(q+1)^{-n}\,\varphi^{(nr)}(z)}{(nr)!(n!)^{-1}}\,t^{n(q+1)}=\sum_{n=0}^{\infty}\frac{\bigg(a^{\frac{1}{r}}(q+1)^{-\frac{1}{r}}t^{\frac{q+1}{r}}\bigg)^{nr}\varphi^{(nr)}(z)}{(nr)!}\\
=\frac{1}{r}\sum_{j=0}^{r-1}\varphi\bigg(z+ e^{\frac{2\pi j}{r}i}a^{\frac{1}{r}}(q+1)^{-\frac{1}{r}}t^{\frac{q+1}{r}}\bigg).
\end{multline*}

 Let $f(s,z):=(\Bo_m\widehat{u})(s,z)$, then by using $m$-moment Laplace transform in a nonsingular direction $d$ we get 
\begin{multline*} (T_{m,d}f)(t,z)=\int_{e^{id}\RR_+}e_m\Big(\frac{s}{t}\Big)\*f(s,z)\frac{ds}{s}\\
 =\frac{1}{r}\int_{e^{id}\RR_+}e_m\Big(\frac{s}{t}\Big)\*\*\sum_{j=0}^{r-1}\varphi\bigg(z+ e^{\frac{2\pi j}{r}i}a^{\frac{1}{r}}(q+1)^{-\frac{1}{r}}s^{\frac{q+1}{r}}\bigg)\frac{ds}{s}.
 \end{multline*}

Thus,  the unique formal solution $\widehat{u}(t,z)$ of this Cauchy problem is $\frac{q+1}{r-1}$-summable in the direction $d$ and for every $\varepsilon>0$ there exists $\tilde r>0$ such that
its $\frac{q+1}{r-1}$-sum $u\in\OO(S_{d}(\pi\frac{r-1}{q+1}-{\varepsilon},\tilde r)\times D)$ is given by 
$$
u(t,z)=u^{d}(t,z)\\=\frac{1}{r}\*\int_{e^{i d}\RR_+}e_m\Big(\frac{s}{t}\Big)\*\*\sum_{j=0}^{r-1}\varphi\bigg(z+ e^{\frac{2\pi j}{r}i}a^{\frac{1}{r}}(q+1)^{-\frac{1}{r}}s^{\frac{q+1}{r}}\bigg)\frac{ds}{s}.
$$

Let $\theta_z:=\arg (z_0-z)$, $\delta_l:=\frac{r\theta_z-2\pi l-\arg a}{q+1}$, where $l=0,1,\ldots,r-1$. 
Then  $\mathcal{L}_{\delta_l}$    are Stokes lines for $\widehat{u}$.
For every sufficiently small $\varepsilon>0$ there exists $r>0$ such that for every fixed $z\in D_r$ the jump is given by 
\begin{multline*}
J_{\mathcal{L}_{\delta_l}}\widehat u(t,z)= u^{\delta_l+\varepsilon}(t,z) - u^{\delta_l-\varepsilon}(t,z)=F_z(s)\bigg[\frac{e_m(s/t)}{s}\bigg]=\big[\Bo_m\widehat{u}(s,z)\big]\bigg[\frac{e_m(s/t)}{s}\bigg]\\
=\biggl[\sum_{j=0}^{r-1}\varphi\bigg(z+ e^{\frac{2\pi j}{r}i}a^{\frac{1}{r}}(q+1)^{-\frac{1}{r}}s^{\frac{q+1}{r}}\bigg)\biggr]_{\delta_l}\biggl[\frac{e_m(s/t)}{rs}\biggr]\\
=\biggl[\varphi\big(z+e^{\frac{2\pi l}{r}i}a^{\frac{1}{r}}(q+1)^{-\frac{1}{r}}s^{\frac{q+1}{r}}\big)\biggr]_{\delta_l}\biggl[\frac{e_m(s/t)}{rs}\biggr].
\end{multline*}

The last equality arising from the fact that in this case all singular points appear in the function $s\mapsto\varphi\big(z+e^{\frac{2\pi l}{r}i}a^{\frac{1}{r}}(q+1)^{-\frac{1}{r}}s^{\frac{q+1}{r}}\big)$.

Observe that from \cite[Theorem 32]{B2} one can derive the function
\begin{multline*}
 e_m(u)=T^{-}_{m_1,d}\,\bigg(e_{m_2}\Big(\frac{1}{z}\Big)\bigg)\Big(\frac{1}{u}\Big)=-\frac{1}{2\pi i}\int_{\gamma(d)}E_{m_1}\Big(\frac{1}{uz}\Big)e_{m_2}\Big(\frac{1}{z}\Big)\frac{dz}{z}\\
 \stackrel{\mathrm{by\,\,\, Example} \ref{ex:functions}}{=}-\frac{1}{2\pi i}\int_{\gamma(d)}\sum_{v=0}^{\infty}\frac{(uz)^{-v}}{m_1(v)}\,\,\frac{q+1}{r}\Big(\frac{1}{z}\Big)^\frac{q+1}{r}e^{-\big(\frac{1}{z}\big)^\frac{q+1}{r}}\frac{dz}{z}\\
 =-\frac{1}{2\pi i}\int_{\gamma(d)}\mathbf{E}_{\frac{1}{q+1}}\Big(\frac{1}{uz}\Big)\,\,\frac{q+1}{r}\Big(\frac{1}{z}\Big)^\frac{q+1}{r}e^{-\big(\frac{1}{z}\big)^\frac{q+1}{r}}\frac{dz}{z}.
\end{multline*}

\item Let now $r>1,p\neq0$ and $\varphi(z)\in\OO^{\frac{r}{p-1+r}}\Bigl(\widetilde{\CC\setminus\{z_0\}}\Bigr)$ for some $z_0\in\CC\setminus\{0\}$. Then by Proposition \ref{propo:1} function $m(n):=m_2(n)\cdot\underbrace {m_1(n)\cdot\ldots\cdot m_1(n)}_{(p-1)-times}=\Gamma\big(1+\frac{nr}{q+1}\big)\bigg(\Gamma\big(1+\frac{n}{q+1}\big)\bigg)^{p-1}$  is a moment function of order $\frac{p-1+r}{q+1}$ corresponding to a kernel function $e_{m}(z)$ of order $\frac{q+1}{p-1+r}$.

We start with applying the $m$-moment Borel transform to $\widehat{u}(t,z)$. 
\begin{multline*}
(\Bo_m\widehat{u})(t,z)=\sum_{n=0}^{\infty}\frac{a^n (n!)^{p-1}\,(q+1)^{n(p-1)}\,\varphi^{(nr)}(z)}{m(n(q+1))}\,t^{n(q+1)}\\
=\sum_{n=0}^{\infty}\frac{a^n (n!)^{p-1}\,(q+1)^{n(p-1)}\,\varphi^{(nr)}(z)}{\Gamma\big(1+\frac{n(q+1)r}{q+1}\big)\bigg(\Gamma\big(1+\frac{n(q+1)}{q+1}\big)\bigg)^{p-1}}\,t^{n(q+1)}=\sum_{n=0}^{\infty}\frac{a^n (n!)^{p-1}\,(q+1)^{n(p-1)}\,\varphi^{(nr)}(z)}{\Gamma(1+nr)\bigg(\Gamma(1+n)\bigg)^{p-1}}\,t^{n(q+1)}\\
=\sum_{n=0}^{\infty}\frac{a^n (n!)^{p-1}\,(q+1)^{n(p-1)}\,\varphi^{(nr)}(z)}{(nr)!(n!)^{p-1}}\,t^{n(q+1)}=\sum_{n=0}^{\infty}\frac{\bigg(a^{\frac{1}{r}}(q+1)^{\frac{p-1}{r}}t^{\frac{q+1}{r}}\bigg)^{nr}\varphi^{(nr)}(z)}{(nr)!}\\
=\frac{1}{r}\sum_{j=0}^{r-1}\varphi\bigg(z+ e^{\frac{2\pi j}{r}i}a^{\frac{1}{r}}(q+1)^{\frac{p-1}{r}}t^{\frac{q+1}{r}}\bigg).
\end{multline*}

 Let $f(s,z):=(\Bo_m\widehat{u})(s,z)$, then by using $m$-moment Laplace transform in a nonsingular direction $d$ we get 
\begin{multline*} (T_{m,d}f)(t,z)=\int_{e^{id}\RR_+}e_m(s/t)\*f(s,z)\frac{ds}{s}\\
 =\frac{1}{r}\int_{e^{id}\RR_+}e_m(s/t)\*\*\sum_{j=0}^{r-1}\varphi\bigg(z+ e^{\frac{2\pi j}{r}i}a^{\frac{1}{r}}(q+1)^{\frac{p-1}{r}}s^{\frac{q+1}{r}}\bigg)\frac{ds}{s}.
 \end{multline*}

Thus,  the unique formal solution $\widehat{u}(t,z)$ of this Cauchy problem is $\frac{q+1}{p-1+r}$-summable in the direction $d$ and for every $\varepsilon>0$ there exists $\tilde r>0$ such that
its $\frac{q+1}{p-1+r}$-sum $u\in\OO(S_{d}(\pi\frac{p-1+r}{q+1}-{\varepsilon},\tilde r)\times D)$ is given by 
$$
u(t,z)=u^{d}(t,z)\\=\frac{1}{r}\*\int_{e^{i d}\RR_+}e_m(s/t)\*\*\sum_{j=0}^{r-1}\varphi\bigg(z+ e^{\frac{2\pi j}{r}i}a^{\frac{1}{r}}(q+1)^{\frac{p-1}{r}}s^{\frac{q+1}{r}}\bigg)\frac{ds}{s}.
$$

Let $\theta_z:=\arg (z_0-z)$, $\delta_l:=\frac{r\theta_z-2\pi l-\arg a}{q+1}$, where $l=0,1,\ldots,r-1$. 
Then  $\mathcal{L}_{\delta_l}$ are Stokes lines for $\widehat{u}$.
For every sufficiently small $\varepsilon>0$ there exists $\tilde r>0$ such that for every fixed $z\in D_{\tilde r}$ the jump is given by 
\begin{multline*}
J_{\mathcal{L}_{\delta_l}}\widehat u(t,z)= u^{\delta_l+\varepsilon}(t,z) - u^{\delta_l-\varepsilon}(t,z)=F_z(s)\bigg[\frac{e_m(s/t)}{s}\bigg]=\big[\Bo_m\widehat{u}(s,z)\big]\bigg[\frac{e_m(s/t)}{s}\bigg]\\
=\biggl[\sum_{j=0}^{r-1}\varphi\bigg(z+ e^{\frac{2\pi j}{r}i}a^{\frac{1}{r}}(q+1)^{\frac{p-1}{r}}s^{\frac{q+1}{r}}\bigg)\biggr]_{\delta_l}\biggl[\frac{e_m(s/t)}{rs}\biggr]\\
=\biggl[\varphi\big(z+e^{\frac{2\pi l}{r}i}a^{\frac{1}{r}}(q+1)^{\frac{p-1}{r}}s^{\frac{q+1}{r}}\big)\biggr]_{\delta_l}\biggl[\frac{e_m(s/t)}{rs}\biggr]
\end{multline*}

The last equality arising from the fact that in this case all singular points appear in the function $s\mapsto\varphi\big(z+e^{\frac{2\pi l}{r}i}a^{\frac{1}{r}}(q+1)^{\frac{p-1}{r}}s^{\frac{q+1}{r}}\big)$.

Notice that from \cite[Theorem 31]{B2} one can derive the function $e_m(z)$ by induction on $p$ ie.
\begin{itemize}
\item Set $p=1$, then $m^{\mathbf{0}}:=m(n)= m_2(n)$ and 
$$ e_ {m^{\mathbf{0}}}(z)= e_{m_2}(z)\stackrel{\mathrm{by\,\,\, Example}\, \ref{ex:functions}}{=}\frac{q+1}{r}\,\,z^{\frac{q+1}{r}}\,e^{-z^{\frac{q+1}{r}}}.$$
\item Set $p=2$, then $m^{\mathbf{1}}(n):=m(n)=m^{\mathbf{0}}(n)m_1(n)$ and 
\begin{multline*}
 e_{m^{\mathbf{1}}}(z)=T_{m_1,d}\,\bigg(e_{m^{\mathbf{0}}}\Big(\frac{1}{u}\Big)\bigg)\Big(\frac{1}{z}\Big) =\int_{e^{id}\RR_+} e_{m_1}(uz)\,\,e_{m^{\mathbf{0}}}\Big(\frac{1}{u}\Big)\frac{du}{u}\\
 \stackrel{\mathrm{by\,\,\, Example}\, \ref{ex:functions}}{=}\int_{e^{id}\RR_+}(q+1)(uz)^{q+1}e^{-(uz)^{q+1}}\frac{q+1}{r}\,\Big(\frac{1}{u}\Big)^{\frac{q+1}{r}}\,e^{-\big(\frac{1}{u}\big)^{\frac{q+1}{r}}}\frac{du}{u}.
\end{multline*}

\item Set $p=3$, then $m^{\mathbf{2}}(n):=m(n)=m^{\mathbf{1}}(n)m_1(n)$ and 
\begin{multline*}
 e_{m^{\mathbf{2}}}(z)=T_{m_1,d}\,\bigg(e_{m^{\mathbf{1}}}\Big(\frac{1}{u}\Big)\bigg)\Big(\frac{1}{z}\Big) =\int_{e^{id}\RR_+} e_{m_1}(uz)\,\,e_{m^{\mathbf{1}}}\bigg(\frac{1}{u}\bigg)\frac{du}{u}\\
 \stackrel{\mathrm{by\,\,\, Example}\, \ref{ex:functions}}{=}\int_{e^{id}\RR_+}(q+1)(uz)^{q+1}e^{-(uz)^{q+1}}e_{m^{\mathbf{1}}}\bigg(\frac{1}{u}\bigg)\frac{du}{u}.
 \end{multline*}
 
 \item Hence, in general if $p\geq2$, then $m^{\mathbf{p-1}}(n):=m(n)=m^{\mathbf{p-2}}(n)m_1(n)$ and
 \begin{multline*}
  e_{m^{\mathbf{p-1}}}(z)=T_{m_1,d}\,\bigg(e_{m^{\mathbf{p-2}}}\Big(\frac{1}{u}\Big)\bigg)\Big(\frac{1}{z}\Big) =\int_{e^{id}\RR_+} e_{m_1}(uz)\,\,e_{m^{\mathbf{p-2}}}\bigg(\frac{1}{u}\bigg)\frac{du}{u}\\
 \stackrel{\mathrm{by\,\,\, Example}\, \ref{ex:functions}}{=}\int_{e^{id}\RR_+}(q+1)(uz)^{q+1}e^{-(uz)^{q+1}}e_{m^{\mathbf{p-2}}}\bigg(\frac{1}{u}\bigg)\frac{du}{u}.
 \end{multline*}

\end{itemize}
\end{enumerate}

\bigskip

Finally, we summarize all derivations and write our main results.

\begin{Tw}
\label{tw1}
Suppose that $d\neq\frac{2k\pi-\arg a}{q+1}$ for $k\in\ZZ$. Assume also that $\hat{u}(t,z)$ is a unique formal solution of the Cauchy problem of the equation (\ref{eq}) with $r=0$ and $p\geq2$.
Then $\hat{u}(t,z)$ is $\frac{q+1}{p-1}$-summable in the direction $d$ and for every $\varepsilon>0$ there exists $\tilde r>0$ such that its $\frac{q+1}{p-1}$-sum $u\in\OO(S_{d}(\frac{\pi(p-1)}{q+1}-{\varepsilon},\tilde r)\times D)$ is given by 
$$
u(t,z)=u^{d}(t,z)=(T_{m,d}\,\Bo_m\widehat{u})(t,z)=\int_{e^{id}\RR_+}e_m\Big(\frac{s}{t}\Big)\*\*\frac{\varphi(z)}{1-a\,(q+1)^{p-1}\,s^{q+1}}\frac{ds}{s}.
$$
Moreover, if $d_k:=\frac{2k\pi-\arg{a}}{q+1}$ for $k=0,\ldots,q$, then  $\mathcal{L}_{d_k}$ are Stokes lines for $\widehat{u}$ and the jump is given by 
$$
J_{\mathcal{L}_{d_k}}\widehat u(t,z)=2\pi i \frac{\varphi(z)}{q+1} \frac{e_m\Big(\big(\frac{1}{a}\big)^{\frac{1}{q+1}}(q+1)^{-\frac{p-1}{q+1}}\,e^{-\frac{2\pi ik}{q+1}}t^{-1}\Big)}{\big(\frac{1}{a}\big)^{\frac{1}{q+1}}(q+1)^{-\frac{p-1}{q+1}}\,e^{-\frac{2\pi ik}{q+1}}}.
$$
\end{Tw}

\begin{Tw}
\label{tw2}
Let $\hat{u}(t,z)$ be a unique formal solution of the Cauchy problem of the equation (\ref{eq}) with $p\geq0$, $r>1$ and $\varphi(z)\in\OO^{\frac{r}{p-1+r}}\Bigl(\widetilde{\CC\setminus\{z_0\}}\Bigr)$ for some $z_0\in\CC\setminus\{0\}$.
Then $\widehat{u}(t,z)$  is $\frac{q+1}{p-1+r}$-summable in the direction $d$ and for every $\varepsilon>0$ there exists $\tilde r>0$ such that
its $\frac{q+1}{p-1+r}$-sum $u\in\OO(S_{d}(\pi\frac{p-1+r}{q+1}-{\varepsilon},\tilde r)\times D)$ is given by 
$$
u(t,z)=u^{d}(t,z)=(T_{m,d}\,\Bo_m\widehat{u})(t,z)=\frac{1}{r}\*\int_{e^{i d}\RR_+}e_m(s/t)\*\*\sum_{j=0}^{r-1}\varphi\bigg(z+ e^{\frac{2\pi j}{r}i}a^{\frac{1}{r}}(q+1)^{\frac{p-1}{r}}s^{\frac{q+1}{r}}\bigg)\frac{ds}{s}.
$$
Moreover, if $\theta_z:=\arg (z_0-z)$, $\delta_l:=\frac{r\theta_z-2\pi l-\arg a}{q+1}$, where $l=0,1,\ldots,r-1$, then $\mathcal{L}_{\delta_l}$ are Stokes lines for $\widehat{u}(t,z)$ and for every sufficiently small $\varepsilon>0$ there exists $\tilde r>0$ such that for every fixed $z\in D_{\tilde r}$ the jump is given by 
$$
J_{\mathcal{L}_{\delta_l}}\widehat u(t,z)= \biggl[\varphi\big(z+e^{\frac{2\pi l}{r}i}a^{\frac{1}{r}}(q+1)^{\frac{p-1}{r}}s^{\frac{q+1}{r}}\big)\biggr]_{\delta_l}\biggl[\frac{e_m(s/t)}{rs}\biggr].
$$
\end{Tw}
\end{section}

\bibliographystyle{siam}
\bibliography{summa}

\end{document}